\documentclass[english]{article}       
\usepackage{graphicx}
\usepackage{latexsym}
\usepackage{amssymb}
\usepackage[mathscr]{eucal}

\newcommand{\Ep}{\hfill\hfill{\qed}\end{proof}}
\newcommand{\Bp}{\begin{proof}}

\newcommand{\Div}{\mbox{\rm div}\,}

\newcommand{\Int}[2]{{\displaystyle \int_{ #1}^{ #2}}}

\newcommand{\Sum}[2]{{\displaystyle \sum_{#1}^{#2}}}

\newcommand{\Frac}[2]{\displaystyle{\frac{\displaystyle{#1}}{\displaystyle{#2}}}}

\newcommand{\beea}{\begin{eqnarray}}

\newcommand{\eeea}{\end{eqnarray}}

\newcommand{\bfe}{{\mbox{\boldmath $e$}} }

\newcommand{\0}{{\mbox{\boldmath $0$}} }

\newcommand{\BF}{\begin{footnotesize}}

\newcommand{\EF}{\end{footnotesize}}

\newcommand{\bi}{\begin{itemize}}

\newcommand{\ei}{\end{itemize}}

\newcommand{\ed}{\end{document}}

\newcommand{\be}{\begin{equation}}

\newcommand{\ba}{\begin{array}}

\newcommand{\ea}{\end{array}}

\newcommand{\ee}{\end{equation}}

\newcommand{\eeq}[1]{\label{eq:#1}\end{equation}}

\newcommand{\real}{{\mathbb R}}

\newcommand{\bfx}{\mbox{\boldmath $x$}}

\newcommand{\bfy}{\mbox{\boldmath $y$}}

\newcommand{\bfv}{{\mbox{\boldmath $v$}} }

\newcommand{\bfu}{{\mbox{\boldmath $u$}} }

\newcommand{\bfw}{{\mbox{\boldmath $w$}} }

\newcommand{\bff}{{\mbox{\boldmath $f$}} }

\newcommand{\bfa}{{\mbox{\boldmath $a$}} }

\newcommand{\bfA}{{\mbox{\boldmath $A$}} }

\newcommand{\bfB}{{\mbox{\boldmath $B$}} }

\newcommand{\bfV}{{\mbox{\boldmath $V$}} }

\newcommand{\bfb}{{\mbox{\boldmath $b$}} }

\newcommand{\bfn}{{\mbox{\boldmath $n$}} }

\def\Bbb R{\real}

\def\tilde{\widetilde}

\newcommand{\bfcalr}{\mbox{\boldmath ${\cal R}$}}

\newcommand{\ED}{\end{description}}

\def\tag{\renewcommand{\theequation}}

\newcommand{\Footnote}{~\footnote}

\newcommand{\Br}{\begin{rem}\begin{rm}}

\newcommand{\Er}{\end{rm}\end{remark}}

\newtheorem{lemm}{Lemma}[section]

\newtheorem{theo}{Theorem}[section]

\newtheorem{rem}{Remark}[section]

\newtheorem{coro}{Corollary}[section]

\newtheorem{exe}{\footnotesize{Exercise}}[section]

\newcommand{\Be}{\begin{exe}\begin{footnotesize}\begin{rm}}

\newcommand{\EE}[1]{\end{rm}\end{footnotesize}\label{exe:#1}\end{exe}}

\newcommand{\Bt}{\begin{theo}\begin{sl}}

\newcommand{\Et}{\end{sl}\end{theorem}}

\newcommand{\Bl}{\begin{lemm}\begin{sl}}

\newcommand{\El}{\end{sl}\end{lemma}}

\newcommand{\eqref}[1]{{\rm (\ref{eq:#1})}}

\newcommand{\Bc}{\begin{coro}\begin{sl}}

\newcommand{\Ec}{\end{sl}\end{coro}}

\newcommand{\ET}[1]{\end{sl}\label{theo:#1}\end{theo}}

\newcommand{\EL}[1]{\end{sl}\label{lemm:#1}\end{lemm}}

\newcommand{\ER}[1]{\end{rm}\label{rem:#1}\end{rem}}

\newcommand{\EC}[1]{\end{sl}\label{coro:#1}\end{coro}}

\begin{document}
\title
{On the Leray-Hopf Extension Condition for the Steady-State Navier--Stokes Problem in Multiply-Connected Bounded Domains}\author{Giovanni P. Galdi}


\maketitle

\begin{abstract}
Employing the approach of A. Takeshita [Pacific J. Math., {\bf 157} (1993), 151--158], we give an elementary proof of the invalidity of the Leray-Hopf Extension Condition for certain multiply connected bounded domains of $\mathbb R^n$, $n=2,3$, whenever the flow through the different components of the boundary is non-zero. Our proof is alternative to  and, to an extent, more direct than the recent one proposed by  J.G. Heywood [J. Math. Fluid Mech. {\bf 13} (2011), 449--457].
\end{abstract}

\section{Introduction}
\label{intro}
Let $\Omega$ be a bounded domain of $\real^n$, $n=2,3$. As is well known, the Leray-Hopf Extension Condition is related to the solvability of the following Navier--Stokes equations
\be \left.\ba{ll}%
\smallskip
\nu\Delta\bfv=\bfv\cdot\nabla\bfv+\nabla p+\bff\\
\Div\bfv=0\ea\right\}\ \ \mbox{in $\Omega$}\,,
\eeq{1}
under prescribed non-homogeneous boundary conditions
\be
\bfv=\bfv_*\ \ \ \mbox{at $\partial\Omega$}.
\eeq{2}
Here, as customary, $\bfv$, $p$ and $\nu>0$ denote velocity and pressure fields, and kinematic viscosity of the liquid, respectively, while $\bff$ is representative of a body force possibly acting on it. Moreover, $\bfv_*$ is a given distribution of velocity at the boundary $\partial\Omega$, which, by \eqref{1}$_1$ and the Gauss theorem must satisfy the compatibility condition
\be  
\Sum{k=1}{N}\int_{\Gamma_k}\bfv_*\cdot\bfn:=
\Sum{k=1}{N}\Phi_k=0\,,
\eeq{3}
where $\Gamma_k$, $k=1,\cdots,N$, are the connected components of $\partial\Omega$, and $\bfn$ is its  unit outer normal. From the physical viewpoint, $\Phi_k$ is the (mass) flow-rate through the portion $\Gamma_k$ of the boundary. To fix the ideas, we assume that  $\Gamma_i$, $i=1,2,\cdots,N-1$, are all surrounded by $\Gamma_N$ and lie outside of each other.  

The existence of a (weak, in principle) solution to the problem \eqref{1}--\eqref{3} is readily established (e.g., by Galerkin method or by Leray-Schauder theory),  provided we are able to show (formally, at least) that the velocity field of the searched solution satisfies the a priori bound
\be
\|\nabla\bfv\|_2\le C\,,
\eeq{4}
where $\|\cdot\|_2$ is the $L^2(\Omega)$-norm,\footnote{We employ standard notation for  Lebesgue, Sobolev and trace spaces; see e.g. \cite{AF}.} and $C$, here and in the following, denotes a constant depending at most on $\Omega$, $\bff$, $\bfv_*$ and $\nu$; see \cite[Chapter IX]{Ga} for details.

One way of attempting to prove \eqref{4} is to extend  the boundary data $\bfv_*$ to $\Omega$ by a solenoidal function $\bfV$, and introduce the new velocity field $\bfu:= \bfV-\bfv$. Clearly, $\bfv$ satisfies a bound of the type \eqref{4} if and only if $\bfu$ does. Now, writing \eqref{1} in terms of $\bfu$, dot-multiplying both sides of the resulting equation by $\bfu$, integrating by parts over $\Omega$ and using the fact that $\bfu$ is solenoidal and that vanishes at $\partial\Omega$, we formally show the following relation:
\be 
\nu\|\nabla\bfu\|_2^2=-(\bfu\cdot\nabla\bfV,\bfu)-
(\bfV\cdot\nabla\bfV,\bfu)-\nu(\nabla\bfV,\nabla\bfu)+(\bff,\bfu)\,,
\eeq{5}
where we have adopted the standard notation\Footnote{Summation convention over repeated indeces applies.}
$$
(\bfa,\bfb)=\Int\Omega{} a_ib_i\,,\ \ \bfa,\bfb\in \real^n\,;\ \ 
(\bfA,\bfB)=\Int\Omega{} A_{ij}B_{ij}\,,\ \ \bfA,\bfB\in \real^{n^2}\,.
$$
Thus, assuming $\bfV\in W^{1,2}(\Omega)$,\footnote{This condition on $\bfV$ is  certainly satisfied if $\bfv_*\in W^{1/2,2}(\partial\Omega)$ and $\Omega$ is Lipschitz.\label{foot:2}} and using  Cauchy--Schwartz and classical embedding inequalities, from \eqref{5} we show
\be
\nu\|\nabla\bfu\|_2^2\le -(\bfu\cdot\nabla\bfV,\bfu)+C\,.
\eeq{6}
Since both terms in \eqref{6} are quadratic in $\bfu$, from this relation it is not clear how to get a bound on $\nabla\bfu$ of the type \eqref{4}, unless we make the obvious assumption that the viscosity $\nu$ is ``sufficiently large'' compared to the magnitude of $\nabla\bfV$, or equivalently,  of $\bfv_*$ in suitable trace norm. However, such a restriction can be avoided whenever $\bfv_*$ obeys the {\em Leray--Hopf Extension Condition}  \cite[p. 38]{Leray}, \cite[p. 772]{Hopf}, namely,  {\em for any $\varepsilon>0$, there exists a solenoidal extension, $\bfV_\varepsilon\in W^{1,2}(\Omega)$,\Footnote{See footnote \ref{foot:2}.} of $\bfv_*$ such that
\tag{EC}
\be
-(\bfu\cdot\nabla\bfV_\varepsilon,\bfu)\le \varepsilon\|\nabla\bfu\|_2^2\,,
\label{EC}
\ee
\renewcommand{\theequation}{\arabic{equation}}\setcounter{equation}{6}for all solenoidal vector functions} $\bfu\in W^{1,2}_0(\Omega)$.
\footnote{Notice that (\ref{EC}) is {\em weaker} than the so-called ``Leray  Inequality'', the latter consisting in replacing the left-hand side of (\ref{EC}) with $|(\bfu\cdot\nabla\bfV_\varepsilon,\bfu)|$. The validity of Leray's Inequality is originally studied, and disproved under certain conditions, in \cite{Tak} and, more recently, in \cite{FKY}.} 
It is then obvious that the validity of (\ref{EC}) along with \eqref{6} furnishes the desired uniform bound for $\bfu$, without imposing any restriction on the magnitude of $\nu>0$.

The validity of (\ref{EC}) has been investigated by many authors, beginning with the cited pioneering works of J. Leray and E. Hopf, who showed that (\ref{EC}) certainly holds provided $\bfv_*$ satisfies a condition {\em stronger} than \eqref{3}, namely, that $\Phi_k=0$, for each $k=1,\cdots,N$. More  recently, a proof of (\ref{EC}) under the {\em general} assumption \eqref{3} was given  in the two-dimensional case  by L.I. Sazonov \cite{Sazonov}, and, independently, by H. Fujita \cite{Fujita}, provided, however $\Omega$, $\bfv_*$, and $\bfu$  satisfy suitable {\em symmetry} hypotheses; see also \cite{Morimoto,FM}. As a result, existence  to problem \eqref{1}--\eqref{3} follows on condition that also $\bff$ is prescribed in an appropriate class of symmetric functions.  The method of Fujita was successively extended by V.V. Pukhnachev \cite{P} to cover the three-dimensional case, again under appropriate symmetry assumptions. 
      
The fact that (\ref{EC}) may not be true unless some restrictions are imposed,  was already clear after the work of A. Takeshita \cite[Section 3]{Tak} and the present author \cite[pp. 22--23]{Ga0}, where it was shown that even in the simplest case when $\Omega$ is an annulus $\mathscr A$, (\ref{EC}) fails in general. More precisely, denoting by $\Gamma_2$ and $\Gamma_1$ the outer and inner concentric circles bounding $\mathscr A$,  one proves that (\ref{EC}) cannot hold at least  when the flow-rate, $\Phi:=\Phi_2=-\Phi_1$ through $\Gamma_2$ is {\em strictly negative} (inflow condition).  A similar result remains valid also in the case when $\Omega$ is a spherical shell, as stated in \cite[p. 157]{Tak} and clearly worked out in \cite{FKY}.

The counterexamples mentioned above require $\Phi<0$. The case $\Phi>0$ (outflow condition) presented, presumably, more difficulty and, as a result, the question of whether (\ref{EC}) holds under the latter assumption on the flow-rate remained apparently open for several years.\Footnote{In this regard, see \cite[Section 8]{He}.} Quite recently, in \cite{He}, J.G. Heywood finally provided very interesting ideas on how to show the invalidity of (\ref{EC}) also for the case $\Phi>0$. This is achieved by using appropriate functions $\bfu$ in (\ref{EC}), that he names ``U-tube test functions.''    

Objective of this note is to give a direct and elementary proof of the invalidity of (\ref{EC}) when $\Omega$ is an annulus (see Section 2) or a spherical shell (see Section 3), and $\Phi>0$.  
Our proof uses Takeshita's approach --which allows us to replace in (\ref{EC}) the extension $\bfV_\varepsilon$ with its integral average over all possible rotations-- in conjunction with an appropriate choice of the function $\bfu$. 

It should be emphasized that once the result is established for these special domains,  the  invalidity of (\ref{EC}) can be extended to more general domains, even multiply connected, whenever for each ``interior'' connected component $\Gamma_i$ of $\partial\Omega$, $i=1,2,\ldots,N-1$, there is a circumference (spherical surface) completely contained in $\Omega$, and that surrounds only $\Gamma_i$. Actually, combining the results of \cite{FKY} with ours, no restrictions need to be imposed on the sign of $\Phi_k$, provided, of course, \eqref{3} is satisfied. This generalization can be obtained by following exactly the same argument of \cite[Corollary 1]{FKY}, and it is stated in Theorem 1  in Section 4.

We wish to end this introductory section with a final observation. A different way of proving the  a priori estimate \eqref{4}, again suggested by J. Leray \cite[p. 28 and {\em ff}\,]{Leray}, is to use a contradiction argument. By this argument one shows that \eqref{4} is true (and so existence to \eqref{1}--\eqref{2} is proved under the general condition \eqref{3}) provided the following requirements on the pair $(\bfw,\pi)$
(in a suitable function class) are incompatible
\be
\ba{cc}\smallskip\left.\ba{ll}\smallskip
\bfw\cdot\nabla\bfw=\nabla\pi\\
\Div\bfw=0\ea\right\}\ \ \mbox{in $\Omega$}\,,\\ \smallskip
\bfw=\0\ \ \mbox{on $\partial\Omega$}\,,\\
-(\bfw\cdot\nabla\bfV,\bfw)=\nu\,. 
\ea
\eeq{7}
Here $(\bfw,\pi)$ are limits (in appropriate topology) of certain normalized sequences of solutions to \eqref{1}--\eqref{2}, while $\bfV$ is a given extension of the boundary data $\bfv_*$. It is then interesting to notice that if \eqref{7} has a solution, then (\ref{EC}) cannot be true. In fact, writing $\bfV_\varepsilon=\bfV+(\bfV_\varepsilon-\bfV)$, from  \eqref{7}$_{1,4}$ we find 
$$
-(\bfw\cdot\nabla\bfV_\varepsilon,\bfw)=\nu+(\bfw\cdot\nabla\bfw,\bfV_\varepsilon-\bfV)=\nu+(\nabla\pi,\bfV_\varepsilon-\bfV)=\nu
$$
where, in the last step, we have used that $\bfV_\varepsilon-\bfV$ is solenoidal and vanishes at $\partial\Omega$. Consequently, admitting (\ref{EC}) would imply 
$\nu\le\varepsilon \|\nabla\bfw\|_2^2$ for all $\varepsilon>0$, namely, $\nu=0$. These considerations suggest that the contradiction argument could be a weaker requirement than the validity of (\ref{EC}), and that it might lead to the proof of the a priori estimate \eqref{4} under more general assumptions than those of symmetry  requested by   the use of (\ref{EC}). This fact was already hinted by C.J. Amick \cite{Amick}, but only recently was it fully confirmed by 
M.V. Korobkov, K. Pileckas, and R. Russo \cite{KPR} who showed  that \eqref{7} are indeed incompatible  when $\Omega$ is a doubly-connected, two-dimensional (Lipschitz) domain, under the sole assumption that the flow-rate satisfies the inflow condition.  
\section{The Case  $\Omega$  an Annulus}
\label{sec:1}
We follow and specialize the approach of \cite{Tak}. 
Let $\Omega:=\{x\in\real^2: R_1<|x|<R_2\}$,  $R_1>0$,  $\Gamma_i:=\{x\in\real^2:|x|=R_i\}$, $i=1,2$. Moreover, set
$$
\Phi:=\Int{\Gamma_2}{}\bfv_*\cdot\bfn=-\Int{\Gamma_1}{}\bfv_*\cdot\bfn\,,
$$
and assume $\Phi>0$. We want to show that the validity of (\ref{EC}) then leads to a contradiction.
For $\bfx\in\Omega$, we put 
\be
\bfy=\bfcalr_{\varphi}\cdot\bfx
\eeq{8}
with $\bfcalr_{\varphi}\in SO(2)$  rotation matrix of angle $\varphi\in [0,2\pi]$, and  define the average of $\bfV_\varepsilon$:
$$
\mathscr A(\bfV_\varepsilon)(\bfy):=\Frac{1}{2\pi}\int_0^{2\pi}
\bfcalr_{\varphi}\cdot\bfV_\varepsilon(\bfcalr_{\varphi}^\top\cdot\bfy)\,d\varphi\,,
$$
where $\bfV_\varepsilon\in W^{1,2}(\Omega)$ is a solenoidal extension of $\bfv_*$ for which (\ref{EC}) is supposed to hold, and ${}^\top$ denotes transpose.
Taking into account the properties of $\bfV_\varepsilon$, and the proper orthogonality of $\bfcalr_\varphi$, one at once shows that
\be
\ba{c}\medskip
\Div \mathscr A(\bfV_\varepsilon)(\bfy)=0\ \ \mbox{$\bfy\in\Omega$}\,,\\
\Int{\Gamma_2}{}\mathscr A(\bfV_\varepsilon)\cdot\bfn=-\Int{\Gamma_1}{}\mathscr A(\bfV_\varepsilon)\cdot\bfn=\Phi\,.
\ea
\eeq{9}
Furthermore, by construction, $\mathscr A(\bfV_\varepsilon)$ is invariant under rotation. Therefore, observing that, denoted by $(r,\theta)$ a system of polar coordinates with the origin at $x=0$, the corresponding base vectors $\{\bfe_r,\bfe_\theta\}$ are both invariant, we infer   
$$
\mathscr A(\bfV_\varepsilon)=v_1(r)\bfe_r+v(r)\bfe_\theta\,.
$$ 
However, $\mathscr A(\bfV_\varepsilon)$ must  satisfy \eqref{9}, so that we conclude  
\be
\mathscr A(\bfV_\varepsilon)=\frac{1}{2\pi}\frac{\Phi}{r}\bfe_r+v(r)\bfe_\theta\,.
\eeq{10}
It is now straightforward to prove that since $\bfV_\varepsilon$ satisfies (\ref{EC}),  also $\mathscr A(\bfV_\varepsilon)$ does. In fact, following \cite{Tak}, by Fubini theorem and \eqref{8},  for all solenoidal $\bfu\in W_0^{1,2}(\Omega)$ we have
\be\ba{rl}\medskip
\Int\Omega{}\bfu\cdot\nabla_y(\mathscr A(\bfV_\varepsilon))\cdot\bfu\, dy=&\!\!\Frac{1}{2\pi}\Int0{2\pi}\Big(\Int\Omega{} \bfu\cdot\nabla_y(\bfcalr_\varphi\cdot\bfV_\varepsilon(\bfcalr_\varphi^\top\cdot\bfy))\cdot\bfu\,dy\Big)d\varphi\\
=&\!\!\Frac{1}{2\pi}\Int0{2\pi}\Big(\Int\Omega{} (\bfcalr_\varphi^\top\cdot\bfu)\cdot\nabla_x\bfV_\varepsilon\cdot (\bfcalr_\varphi^\top\cdot\bfu)\,dx\Big)d\varphi
\ea
\eeq{11}
Clearly, $\bfcalr_\varphi^\top\cdot\bfu\in W_0^{1,2}(\Omega)$ and is solenoidal, and $\|\nabla(\bfcalr_\varphi^\top\cdot\bfu)\|_2=\|\nabla\bfu\|_2$, so that from (\ref{EC})
we deduce
$$
-\Frac{1}{2\pi}\Int0{2\pi}\Big(\Int\Omega{} (\bfcalr_\varphi^\top\cdot\bfu)\cdot\nabla_x\bfV_\varepsilon\cdot (\bfcalr_\varphi^\top\cdot\bfu)\,dx\Big)d\varphi\le\Frac{\varepsilon}{2\pi}\Int0{2\pi} \|\nabla\bfu\|_2^2d\varphi=
\varepsilon \|\nabla\bfu\|_2^2\,.
$$
Combining the latter with \eqref{11} we thus obtain the desired inequality, namely, 
\be
-\Int\Omega{}\bfu\cdot\nabla(\mathscr A(\bfV_\varepsilon))\cdot\bfu\, dy\le\varepsilon\, \|\nabla\bfu\|_2^2\,,
\eeq{12}
for all solenoidal $\bfu\in W_0^{1,2}(\Omega)$. Denoting by $\mathscr D[\mathscr A(\bfV_\varepsilon)]$ the symmetric part of $\nabla \mathscr A(\bfV_\varepsilon)$, we show, on the one hand,
$$
\Int\Omega{}\bfu\cdot\nabla(\mathscr A(\bfV_\varepsilon))\cdot\bfu\, dy= \Int\Omega{}\bfu\cdot\mathscr D[\mathscr A(\bfV_\varepsilon)]\cdot\bfu\, dy 
$$
and, on the other hand, from \eqref{10},
$$
\mathscr D[\mathscr A(\bfV_\varepsilon)]=
\left(\ba{cc}\smallskip -\frac{\Phi}{\pi}\frac{1}{r^2} & v'(r)-\frac{1}{r}v(r)\\
v'(r)-\frac{1}{r}v(r)& \frac{\Phi}{\pi}\frac{1}{r^2} 
\ea\right)\,.
$$
As a result, \eqref{12} becomes
\be
\frac{\Phi}{\pi}\int_\Omega\Frac{u_r^2-u_\theta^2}{r^2}+\int_\Omega (v'(r)-\frac{1}{r}v(r))u_ru_\theta\le \varepsilon\,\|\nabla\bfu\|_2^2\,,
\eeq{13}
for all solenoidal $\bfu\in W_0^{1,2}(\Omega)$, and where $u_r$ and $u_\theta$ denote the polar components of $\bfu$. We now choose $\bfu=(u_r,u_\theta)$ where
\be\ba{ll}\smallskip
u_r=\Frac{m}{r}\Frac{R_2-R_1}{2\pi}\left\{\cos\Big[\Frac{2\pi(r-R_1)}{R_2-R_1}\Big]-1\right\}\cos(m\theta):=U(r)[m\cos(m\theta)]\\
u_\theta=\sin\Big[\Frac{2\pi(r-R_1)}{R_2-R_1}\Big]\sin(m\theta):=W(r)[\sin(m\theta)]\,,\ea
\eeq{14}
with $m$ an integer that will be specified further on. It is obvious that $\bfu\in W^{1,2}_0(\Omega)$, as well as, by taking into account that $(rU)'=-W$,   that $\Div\bfu=0$ in $\Omega$. Moreover, by a direct computation we show (with $\rho=R_2/R_1$)
\be
\int_\Omega\Frac{u_r^2-u_\theta^2}{r^2}=m^2F_1(\rho)-F_2(\rho)\,,
\eeq{15}
where
$$ \ba{ll}\smallskip
F_1(\rho):=\Frac{(\rho-1)^3}{4\pi}\Int01\Frac{1}{[(\rho-1)z+1]^3}\{\cos(2\pi z)-1\}^2dz\,,\\
F_2(\rho):=\pi{(\rho-1)}\Int01\Frac{1}{(\rho-1)z+1}\sin^2(2\pi z)\,dz\,.
\ea
$$
Furthermore, setting $G(r):=(v'(r)-\frac{1}{r}v(r))U(r)W(r)$, we get
\be
\int_\Omega (v'(r)-\frac{1}{r}v(r))u_ru_\theta=m\int_{R_1}^{R_2}rG(r)\,dr\int_0^{2\pi}\sin(m\theta)\cos(m\theta)\,d\theta=0\,.
\eeq{16}
Thus, by fixing $m$ sufficiently large so that
$\kappa:=m^2F_1(\rho)-F_2(\rho)>0$, from \eqref{13}--\eqref{16} we conclude
$$
\kappa\,\Phi\le \varepsilon\, \|\nabla\bfu\|_2^2\,,
$$
which, by the arbitrariness of $\varepsilon>0$, and the assumption $\Phi>0$ furnishes a contradiction. 
\begin{rm} As originally showed by A. Takeshita \cite{Tak}, a similar result also holds if $\Phi<0$. Actually, it is enough to choose in \eqref{13} instead of the field \eqref{14}, the following one
$$
u_r\equiv 0\,,\ \ u_\theta=f(r)
$$
with $f(r)$ any sufficiently smooth function satisfying $f(R_1)=f(R_2)=0$, in which case the left-hand side of \eqref{13} becomes $-{2\Phi}\int_{R_1}^{R_2} rf^2(r)dr$. 
\end{rm}
\section{The case $\Omega$ a Spherical Shell}
In this case $\Omega:=\{x\in\real^3: R_1<|x|<R_2\}$,  $R_1>0$,  $\Gamma_i:=\{x\in\real^3:|x|=R_i\}$, $i=1,2$,
$$
\Phi:=\Int{\Gamma_2}{}\bfv_*\cdot\bfn=-\Int{\Gamma_1}{}\bfv_*\cdot\bfn\,,
$$
and assume $\Phi>0$. Again following the strategy of \cite{Tak}, the proof, in its first part is basically the same as in the case of the two-dimensional annulus. The only change being that the generic rotation matrix is now an element $\bfcalr_{\alpha_1\alpha_2\alpha_3}\in SO(3)$ characterized by the Euler angles $\alpha_i$, $i=1,2,3$. Consequently, \eqref{8} takes the form $\bfy=\bfcalr_{\alpha_1\alpha_2\alpha_3}\cdot\bfx$, and the average $\mathscr A(\bfV_\varepsilon)$ becomes
$$
\mathscr A(\bfV_\varepsilon)(\bfy)=\Frac{1}{8\pi^2}\int_0^{2\pi}\int_0^{2\pi}\int_0^\pi 
\bfcalr_{\alpha_1\alpha_2\alpha_3}\cdot\bfV_\varepsilon(\bfcalr_{_{\alpha_1\alpha_2\alpha_3}}^\top\cdot\bfy)\sin\alpha_3\,d\alpha_1d\alpha_2d\alpha_3\,.
$$
Again, by the properties of $\bfV_\varepsilon$ and the proper orthogonality of the rotation we show that the average satisfies \eqref{9}. Moreover, the invariance of $\mathscr A(\bfV_\varepsilon)$
under the action of  $SO(3)$ along with \eqref{9} implies that
$$
\mathscr A(\bfV_\varepsilon)=\frac{\Phi}{4\pi\,r^2}\,\bfe_r\,,
$$
where  $\{\bfe_r,\bfe_\chi,\bfe_\theta\}$ is the base of a system of spherical coordinates $(r,\chi,\theta)$ with the origin at $x=0$;
see \cite[p. 157]{Tak} for details. Next, proceeding verbatim as in Section 2, we prove that $\mathscr A(\bfV_\varepsilon)$ must satisfy \eqref{12}, which   by taking into account that 
this time 
$$
\{\mathscr D[\mathscr A(\bfV_\varepsilon)]\}_{ij}
=\frac{\Phi}{4\pi r^3}\Big(-3\frac{x_ix_j}{r^2}+\delta_{ij}\Big)\,,
$$
is equivalent to the following
\be
\frac{\Phi}{4\pi}\int_\Omega\frac{2u_r^2-u_{\chi}^2-u_\theta^2}{r^3}\le \varepsilon\,\|\nabla\bfu\|_2^2\,,
\eeq{17}
for all solenoidal $\bfu\in W_0^{1,2}(\Omega)$. In order to show that \eqref{17} leads to a contradiction, we choose
\be
u_r=m\,\tilde{U}(r)\,\cos(m\theta)\sin \chi\,,\ \ u_\chi\equiv 0\,,\ \ u_\theta=\frac{1}{r}W(r)\sin(m\theta)\sin^2\chi\,,
\eeq{18}
where $\tilde{U}=U/r$, and the functions $U$ and $W$ are defined in \eqref{14}. It is easily proved that the vector $\bfu$ with components given in \eqref{18} is solenoidal and is in $W_0^{1,2}(\Omega)$, and therefore can be replaced in \eqref{17}. Since by a direct calculation we show that
$$
\int_\Omega\frac{2u_r^2-u_{\chi}^2-u_\theta^2}{r^3}=m^2G_1(\rho)-G_2(\rho)
$$
where $G_i$, $i=1,2$, are {\em positive} functions of $ \rho=R_2/R_1$, taking $m$ sufficiently large and using the arbitrariness of $\varepsilon$, we show that \eqref{17} is incompatible
with the assumption $\Phi>0$.
\begin{rm} One can show the incompatibility of \eqref{17} also with the alternative assumption $\Phi<0$, by using in \eqref{17} an appropriate  function $\bfu$ different from \eqref{18}. This has been shown in \cite{FKY} by the choice $\bfu=f(r)\sin\chi\bfe_\varphi$, with $f$ sufficiently smooth and satisfying $f(R_1)=f(R_2)=0$, in which case the left-hand side  of \eqref{17} becomes $-\frac{2\Phi}{3}\int_{R_1}^{R_2}r^2f(r)dr$.
\end{rm}
\section{The case $\Omega$ Multiply-Connected}
We now assume that  $\Omega\subset\real^n$, $n=2,3$, is a multiply-connected Lipschitz domain of the type defined in the Introduction. Furthermore, following \cite{FKY}, we suppose that for each  connected component $\Gamma_i$ of $\partial\Omega$, $i=1,2,\ldots,N-1$, there is a circumference (spherical surface) completely contained in $\Omega$, and  surrounding only the component $\Gamma_i$.

Combining the results of the previous two sections with Remarks 1 and 2, and the argument of \cite[Theorem 1 and Corollary 1]{FKY}, we can show the following general result, whose proof follows exactly the same lines of \cite[Corollary 1]{FKY}, and, consequently, will be omitted.
\begin{theo} Let $\Omega\subset\real^n$, $n=2,3$, be a bounded domain satisfying the assumptions mentioned above. Moreover, let $\bfv_*\in W^{1/2,2}(\partial\Omega)$ obey the compatibility condition \eqref{3}. Then the Leray-Hopf Extension Condition holds for $\bfv_*$ (if and) only if $\Phi_k=0$ for all $k=1,2,\ldots,N$.    
\end{theo} 
\medskip\par\noindent
{\bf Acknowledgements.}
Work partially supported by the NSF Grant DMS-1311983



\end{document}